\documentclass[11pt]{article}
\usepackage[pctex32]{graphics}
\usepackage{amssymb}
\usepackage{latexsym}
\usepackage{epsfig} 
\usepackage{multirow} 
\usepackage{pstricks}
\usepackage{float}
\usepackage{epstopdf}
\usepackage{mathrsfs}
\usepackage{soul}

\usepackage{amsmath}
\usepackage{amssymb}
\usepackage{relsize}
\usepackage{multirow}

\usepackage{graphicx}
\usepackage{multicol} 

\usepackage{times}
\usepackage{epstopdf}

\definecolor{Red}{rgb}{1,0.,0.}

\usepackage{fancyhdr}
\usepackage{setspace}

\usepackage{hyperref}
\usepackage{mathtools}
\usepackage{algorithm}
\usepackage[dvipsnames]{xcolor}
\usepackage{algorithmic}
\usepackage{placeins}
\newcommand{\R}{{\mathbb R}}
\newcommand{\N}{{\mathbb N}}

\newcommand{\mS}{{\mathsf S}}

\newcommand{\mP}{{\mathsf P}}
\newcommand{\mT}{{\mathsf T}}
\newcommand{\mA}{{\mathsf A}}

\newcommand{\mR}{{\mathsf R}}
\newcommand{\mI}{{\mathsf I}}
\newcommand{\mC}{{\mathsf C}}
\newcommand{\mD}{{\mathsf D}}

\newcommand{\mW}{{\mathsf W}}

\newcommand{\mV}{{\mathsf V}}
\newcommand{\mU}{{\mathsf U}}

\newcommand{\mK}{{\mathsf K}}
\newcommand{\mL}{{\mathsf L}}
\newcommand{\mQ}{{\mathsf Q}}
\newcommand{\mG}{{\mathsf G}}
\newcommand{\mM}{{\mathsf M}}

\title{Spotlight, priorsketching and Bayesian approximation error paradigms}
\author{D Calvetti \and E Somersalo} 

\date{Case Western Reserve University\\
Department of Mathematics, Applied Mathematics and Statistics\\
Cleveland, OH, USA}

\begin{document}

\maketitle

\begin{abstract}
A way to lower computational cost in large scale inverse problems and problems depending on poorly known model parameters is to replace the detailed model by an approximate one. Inverse problems are typically ill-posed, and the model discrepancy introduced by using  approximate models often shows up in the computed solutions as disturbing artifacts or blurring. In this article, we consider two methods of addressing certain types of modeling errors, the Bayesian approximation error (BAE) method and linear algebraic spotlight inversion to suppress clutter in the computational model by orthogonal projections. Through the process of analyzing the two approaches, we show that they turn out to be closely related but not equivalent, and we highlight a connection to sketching schemes in randomized linear algebra. The similarities between the methods and their successful suppression of most of the clutter effects is elucidated with two computed examples, one addressing of X-ray tomography and the other electrical impedance tomography.   
\end{abstract}
\section{Introduction}

A common challenge in computational inverse problems is to strike a balance between the accuracy and complexity of the computational models for the forward problem. This is especially the case for forward models based on ordinary or partial differential equations when the unknowns to be estimated are distributed parameters.  It is well known that capturing the detailed features of distributed parameters requires dense discretization, leading to a high dimensional representation that may make the inverse solvers very costly, as typically repeated evaluations of the forward model are required. Another source of modeling errors that needs to be addressed in a suitable manner is the dependence of the forward model on poorly known parameters, referred to as nuisance parameters in the statistics literature,  specifying, e.g., the geometric shape of the domain, to which the observed quantity may have a high sensitivity. Nuisance parameters can be included as part of the unknown to be estimated from the data, however this may complicate the inverse problems, e.g., when the data depend linearly on the unknowns of primary interest but non-linearly on the nuisance parameters. 

One way to overcome the curse of dimensionality is to replace the detailed model with computationally less expensive surrogate models of lower dimensionality. Similarly, unknown model parameters may be replaced by fixed values that are believed to represent in some mean sense the reality. It has been shown that resorting to a less accurate forward model introduces a model discrepancy that can be viewed as a modeling error in the data. Due to the inherent ill-posedness of inverse problems, in the inversion process  errors in the data may be enhanced by the ill-conditioning, significantly degrading the computed solutions. Modeling errors are typically quite different from observation errors, often being characterized by nonzero mean and strong correlation among some of the components, which is a reason why neglecting to address their presence may result in the presence of annoying artifacts in the solution. Therefore compensating for the modeling error introduced when replacing the forward model with an approximate models is a very important topic in the computed solution of inverse problems.

The Bayesian framework for inverse problems  \cite{kaipio2005statistical,calvetti2023bayesian}, where all unknowns are modeled as random variables, provides a natural environment to address the modeling errors. Bayesian approximation error (BAE) analysis was originally introduced in \cite{kaipio2005statistical,kaipio2007statistical}, and its effectiveness was illustrated in a wide range of inverse problems; see, e.g., \cite{arridge2006approximation,nissinen2010compensation,nissinen2011reconstruction}\cite{kaipio2013approximate,
koulouri2016compensation,huttunen2007approximation,kaipio2019bayesian}. The idea proposed in \cite{kaipio2007statistical} is to generate a sample of realizations of the  approximation error by means of drawing first a sample from the prior density of the unknown of primary interest and possible nuisance parameters, and computing the differences between the outputs of the accurate and the approximate model. The statistics of the approximation error sample can then be used to adjust the likelihood based on the approximate forward model to account for the approximation error.  This step requires only applications of the forward models, and can be costly due to the computational complexity of the accurate model, however, it can be done offline and only once, akin to the training of machine learning models. 

A different approach for dealing with nuisance parameters in linear inverse problems was recently proposed in \cite{calvetti2025spotlight}, where the observation model is regarded as a linear combination of a term depending linearly on the variable of interest and another with linear dependence on the nuisance parameter. Ideally, the part of no interest, referred to as clutter in the radar literature, can be cleansed by an orthogonal projection, aiming the spotlight on the portion of interest and  reducing the problem to a linear observation model independent of the nuisance parameter. In the context of non-linear models, the projection idea was considered in \cite{jaaskelainen2024projection}. While model reduction by spotlighting is not equivalent to the model reduction compensated by Bayesian approximation error analysis, there is a clear overlap between the two approaches, as pointed out in \cite{calvetti2025spotlight}. The main aim of this article is to investigate in further detail the relationship between the two approaches.

The remainder of the paper is organized as follows. In section 2 we review the Bayesian approximation error approach, and outlining how the formalism can be used in the likelihood to compensate for using an approximate forward operator, and subsequently, we provide an overview of the spotlight inversion. The similarities and differences between the two model reduction approaches are analyzed in section 3, where we first interpret the approximation error as clutter, then interpret the clutter as approximation error. In subsection 3.3 we show that the spotlight inversion can be seen as a draconian version of the Bayesian approximation error adjustment of the likelihood, in which data matching is considered only in the orthogonal complement of the subspace of the approximation error.  In the process  of setting up the machinery for approximation error modeling, we draw some parallels between the Bayesian error modeling and sketching methods that are popular in randomized linear algebra, introducing the concept of {\em priorsketching}.
The Bayesian approximation error approach and spotlight inversion are employed in section 4 in two computed examples. The first computed example is an X-ray tomography problem in two dimensions, where the portion of interest is a subset of the domain representing the target in the imaging domain. The model reduction in this case is obtained by using a fine discretization only for the region of interest and a much coarser one in the remaining part of the domain. The physics informed prior used to estimate the modeling error is based on the logit transformation of a Gaussian prior and tuned according to the a prior beliefs about the solution. The second example considers an electrical impedance tomography (EIT) problem with poorly known geometry of the domain. In this case the clutter is due to the unknown variation of the geometry of the target, and while it cannot be interpreted directly in the spotlight inversion framework as a separate contribution from a nuisance parameter, we show that the priorsketching and orthogonal projection techniques provide a simple and robust way to solve the inverse problem.
Conclusions and future work are the topic of section 5.

\section{Approximation error and spotlight inversion}

We begin this section by reviewing both the Bayesian approximation error paradigm  and the spotlight inversion for clutter reduction. The starting point of the former is the Bayesian formulation of inverse problems, where all unknown parameters are modeled as random variables, while the spotlight inversion is  based on purely linear algebraic considerations. 

\subsection{Reduced model with Bayesian approximation error}
Consider the inverse problem of estimating a variable $x\in\R^N$ from the noisy observation of a related quantity $b\in\R^m$, assuming the forward model
\begin{equation}\label{forward1}
    b = f^*(x) + e,   
\end{equation}
where $f^*:\R^N\to \R^m$ is a known linear or nonlinear function representing the best available approximation of reality, and $e\in\R^m$ represents additive noise. In line with the Bayesian framework, quantities with unknown values are modeled as random variables, the randomness representing the epistemic uncertainty. Using the convention to use uppercase letters to denote random variables and lowercase letters for their realizations, we write the stochastic extension of  (\ref{forward1}) as
\begin{equation}\label{forward2}
B = f^*(X) + E.
\end{equation}
In applications where the evaluation of the forward map $f^*$ is costly,  as is typically the case,  e.g., when it requires the numerical solution of a PDE, it is attractive to replace $f^*$ by a surrogate model of lower computational complexity, obtained by e.g., fixing some unknown parameters of less interest, and/or reducing the dimensionality of the unknown by coarse discretization. For these reasons surrogate models often are referred to as reduced model. Let 
\begin{equation}\label{surrogate1}
 b = f(z) +e', 
\end{equation}
be the reduced model version of (\ref{forward1}), where $z\in\R^n$ is the unknown in the reduced model, $f:\R^n\to\R^m$ is the surrogate forward operator, and $e'$ accounts for measurement errors and for the error induced by the discrepancy between the accurate and the surrogate model, referred to as approximation error. The stochastic extension for the reduced model is 
\begin{equation}\label{surrogate}
    B = f^*(X) + E = f(Z) + \big\{f^*(X) - f(Z)\big\} + E
    =f(Z) + E',
\end{equation}
where the error term $E'$ is of the form
\[
 E' = \big(f^*(X) - f(Z)\big) + E = M+E,
\]
$M$ being the random variable representing the approximation error
\begin{equation}\label{ModelingError}
M = f^*(X) - f(Z). 
\end{equation}
We refer to (\ref{ModelingError}) as Bayesian Approximation Error (BAE).  According to the Bayesian paradigm, $X$ and $Z$ are random variables and their a priori distributions encode the prior belief about the distribution of their values. When replacing the accurate forward operator with the surrogate model, the variable $Z$ is related somehow to the unknown $X$ of the original model, e.g. $Z$ being a projection of $X$ to a lower dimensional subspace. We will assume the existence of a joint prior distribution $\pi_{X,Z}$, in which case the priors for $X$ and $Z$ are the marginals, 
\[
 \pi_X(x) = \int_{\R^n}\pi_{X,Z}(x,z) dz, \quad 
  \pi_Z(z) = \int_{\R^N}\pi_{X,Z}(x,z) dx.
\]
In the Bayesian approximation error paradigm, we define the probability density of the Bayesian approximation error $M$ as the density of the push-forward distribution,
\[
\pi_M = {\mathscr F}_\#\pi_{X,Z}, \quad {\mathscr F}(x,z) = f^*(x) - f(z).
\]
In the case where the distribution of the marginals is not available in closed form, the following  sample-based procedure can be used for the practical estimation of the probability density of the Bayesian approximation error for the inverse problem with the surrogate model (\ref{surrogate}). First, we generate an independently drawn random sample
\begin{equation}\label{sample 1}
\big\{(x^{(1)},z^{(1)}), \ldots, (x^{(L)},z^{(L)})\big\}
\end{equation}
from the joint prior distribution $\pi_{X,Z}$, and use it to generate an associated  approximation error sample,
\begin{equation}\label{sample 2}
 \big\{m^{(1)},\ldots,m^{L)}\big\}, \quad m^{(j)} = { \mathscr F}(x^{(j)},z^{(j)}).
\end{equation}
We then compute the empirical mean and covariance of $M$ based on this sample,
\[
 \mu = \frac 1L\sum_{j=1}^L m^{(j)}, \quad \mC_M = \frac 1L\sum_{j=1}^L(m^{(j)}-\mu)(m^{(j)}-\mu)^\mT,
\]
and we use the Laplace approximation for the distribution of $M$.  Formally,  $M \approx \widetilde M$, where
\begin{equation}\label{laplace}
 \widetilde M\sim{\mathcal N}(\mu,\mC_M).
 \end{equation}
Assuming that $\widetilde M$ is independent of $Z$, we can formulate the surrogate inverse problem as
\begin{equation}\label{bme}
 B = f(Z) + \widetilde M + E.
\end{equation}

If the inverse problem is to be solved in the Bayesian framework, the marginal $\pi_Z$ can be used as a prior for $Z$. If, in addition, the prior for $Z$ is Gaussian,
\[
 Z \sim{\mathcal N}(0,\mC_Z),
\]
the Maximum A Posteriori (MAP) estimate for $Z$ with the Bayesian error model is
\begin{equation}\label{MAP}
z_{\rm MAP} = {\rm argmin}\big\{(b - \mu - f(z))^\mT(\mC_M + \mI_m)^{-1} (b - \mu - f(z)) + z^\mT \mC_Z^{-1} z\big\}.
\end{equation}

Before moving further, a few observations are in order. Originally, when the above approximation error approach was introduced \cite{kaipio2005statistical}, it was referred to as Enhanced Error Model (EEM). One of the tenets of the EEM is the assumed independency of $Z$ and $M$. This assumption is not necessary, and in fact, the interdependency of the samples (\ref{sample 1}) and (\ref{sample 2}) allows the estimation of the cross correlation of $M$ and $Z$, leading to a more detailed error model, see, e.g.,
\cite{calvetti2008hypermodels,kaipio2019bayesian}. Furthermore, the main reason for Laplace approximation (\ref{laplace}) of the probability density of the modeling error is its computational convenience.
A further development of the approximation error model was laid out in
 \cite{calvetti2018iterative}, where
the approximation error model is adaptively updated  as we learn the posterior density of the unknown. 
Since the focus in this article is on connections between the error modeling and spotlight inversion discussed next, we limit our attention to the EEM version of the approximation error.
In the computed examples, we will elaborate further the connections between the original unknown $X$ and the surrogate variable $Z$. 

\subsection{Spotlight inversion}

In  this section, we review the ideas of spotlight inversion in the framework of linear inverse problems. 
Consider a linear inverse problem,
\begin{equation}\label{linsys}
 b = \mA x + e,
\end{equation}
where $\mA\in\R^{m\times N}$ is a known matrix, $x\in\R^N$ is an unknown and $e\in\R^m$ represents additive noise in the observation $b\in\R^m$. In some applications we are interested only in a portion of the vector $x$, which we assume, after possible permutation, to consist of the block of first $n$ entries $x_1\in\R^n$ of $x$, and denote by $x_2\in\R^{N-n}$ the remaining portion of $x$, which we refer to as nuisance parameters. After partitioning the columns of the matrix $\mA$ according to the partitioning of $x$, 
\[
 \mA = \left[\begin{array}{cc} \mA_1 & \mA_2\end{array}\right], \quad x = \left[\begin{array}{c} x_1 \\ x_2 \end{array}
 \right]
\]
where $\mA_1\in\R^{m\times n}$, $\mA_2\in\R^{m\times(N-n)}$, $x_1\in\R^n$, and $x_2\in\R^{N-n}$, we can write (\ref{linsys}) as 
\begin{equation}\label{spot1}
 b = \mA_1 x_1 + \mA_2 x_2 + e.
\end{equation}
Although only $x_1$ is of interest, the nuisance parameter vector $x_2$ contributes to the data $b$ through the {\em clutter term} $\mA_2 x_2$, thus ignoring it typically leads to artifacts or blurring. Recently, we proposed a spotlight inversion algorithm for estimating $x_1$ from the data based on orthogonal projections that eliminate the clutter contribution, see \cite{calvetti2025spotlight}. Below we review the main idea behind the spotlight inversion. 

Assume first that the range ${\mathcal R}(\mA_2)$ of the matrix $\mA_2$ satisfies
\begin{equation}\label{dimension}
 {\rm dim}\left({\mathcal R}(\mA_2)\right)<m,
\end{equation}
and let $\mP$ denote the orthogonal projector onto the range of $\mA_2$,
\[
 \mP: \R^m \to {\mathcal R}(\mA_2),
\]
and by $\mP^\perp$ the associated orthogonal projector onto the orthocomplement of ${\mathcal R}(\mA_2)$,
\begin{equation}\label{P perp}
\mP^\perp = \mI_m - \mP :\R^m\to \big({\mathcal R}(\mA_2)\big)^\perp.
\end{equation}
It is straightforward to verify that if we apply the projector $\mP^\perp$ to both sides of (\ref{spot1}), we obtain the reduced model 
\begin{equation}\label{projected}
 \mP^\perp b = \mP^\perp\mA_1 x_1 +\underbrace{\mP^\perp \mA_2}_{=0} x_2 + \mP^\perp e = \mP^\perp\mA_1 x_1  + \mP^\perp e,
\end{equation}
where the contribution of the nuisance parameter has been completely eliminated. The final step of spotlight inversion entails solving the projected linear inverse problem (\ref{projected}).

While the idea behind the spotlight inversion is very straightforward, there are several potential pitfalls.  The assumption (\ref{dimension}) hold automatically when the columns of the matrix $\mA$ are linearly independent, but it many cases, in particular when the problem is underdetermined, it may not be satisfied. In that case the range of $\mA_2$ is the full space and the projector $\mP^\perp$ maps everything to zero, in which case all information in the data is lost. One way to partially overcome this problem, as suggested in \cite{calvetti2025spotlight}, is to replace $\mP$ with an approximate projector constructed from the first $k$ singular vectors of $\mA_2$.  More specifically,  if
\[
 \mA_2 = \mU \mD \mV^\mT, \quad \mU = \left[\begin{array}{ccc} u_1 & \cdots & u_m\end{array}\right],
\]
we replace $\mP$ in (\ref{P perp}) by an approximation 
\begin{equation}\label{Pk}
 \mP_k = \sum_{j=1}^k u_ju_j^\mT, \quad k<m,
\end{equation}
and define the associated projector  $\mP^\perp$ as $\mP_k^\perp = \mI_m - \mP_k$.  
In the projected equation
\[
 \mP_k^\perp b = \mP_k^\perp \mA_1 x_1 + \big(\mP_k^\perp\mA_2 x_2 + \mP_k^\perp e\big)
\]
there is still a contribution $\mP_k^\perp \mA_2 x_2$ from the projected clutter term, however, in many cases it is possible to select the value of $k$ so that the size of the projected clutter term is at, or below, that of the projected observation error, so that this the latter dominates, and the effect on the inverse solution can be handled by standard regularization techniques.

Another potential complication of the spotlight inversion approach is when the angle between the ranges of $\mA_1$ and $\mA_2$ is very small, in which case the application of the projector $\mP^\perp $ or its approximation $\mP_k^\perp$ may significantly deteriorate the signal-to-noise ratio of the projected problem. Further analysis of these and other potential drawbacks can be found in the cited article. 

In the following section, we explore the similarities and differences between the Bayesian approximation error and the spotlight projection.

\section{ Spotlight on Bayesian approximation error and clutter elimination }

In this section we investigate how to relate  Bayesian approximation error compensation and spotlight inversion. 

\subsection{Spotlighting to curb the approximation error}

Consider the stochastic extension of the reduced model with the approximation error,
\begin{equation}\label{ReducedModel}
 B = f(Z) + M + E.
\end{equation}
For simplicity, we assume that $E$ is zero mean white Gaussian, $E\sim{\mathcal N}(0,\mI_m)$.   Let $\mu$ be the mean of the approximation error $M$, and let $M_c = M- \mu$ be the centered approximation error. 

Consider the covariance matrix of $M$, defined through
\[
 \mG = {\mathbb E}\big(M_c M_c^\mT\big),
\]
and let $\big\{(\lambda_j,u_j)\big\}_{j=1}^m$ be the eigenpairs of $\mG$, $\lambda_1\geq\ldots\geq \lambda_m\geq 0$.
We represent the random variable $M_c$ in the eigenbasis of the covariance matrix as
\[
 M_c = \sum_{j=1}^m \big(u_j^\mT M_c\big) u_j = \sum_{j=1}^m M_j u_j,
\]
where the univariate random variables $M_j$ are uncorrelated since
\[
 {\mathbb E}\big(M_j M_k\big) =
 {\mathbb E}\big(u_j^\mT M_c M_c^\mT u_k\big) = u_j^\mT\mG u_k = \lambda_j \delta_{jk}.
\]
Denoting by $r$ the rank of the matrix $\mG$, we define the random variables $\xi_j$ of unit variance as
\begin{eqnarray*}
 \xi_j & = &  \frac{1}{\sqrt{\lambda_j}} M_j, \quad 1\leq j\leq r \\
       & = & 0 \quad j> r.
\end{eqnarray*}
and obtain a decomposition of the approxmation error of the form
\begin{equation}\label{eig expansion}
 M = \mu + \sum_{j=1}^m \sqrt{\lambda_k} \xi_j u_j \quad \mbox{(a.s.)}.
\end{equation}
This decomposition is the time-invariant analogue of the Karhunen-Lo\'{e}ve decomposition for stochastic processes.

As in the spotlight procedure, consider the projector $\mP_k$ in formula (\ref{Pk}) where  the $u_j$ are the eigenvectors of $\mG$, and apply the complementary projector $\mP_k^\perp$ to (\ref{ReducedModel}).  A comparison of the resulting projected model
\begin{equation}\label{projected eq}
 \mP_k^\perp(B-\mu) = \mP_k^\perp f(Z) + \sum_{j=k+1}^m \sqrt{\lambda_j}\xi_j u_j + \mP_k^\perp E
\end{equation}
with the spotlight analogue shows that the role of $\mA_1 x_1$  is played by the surrogate model $f(z)$ with $z=x_1$. Furthermore, the role of the clutter is taken by the eigenvector expansion of the approximation error, and if the covariance were degenerate and  $\lambda_j = 0$ for $j>k$, in the projected model there would be no contribution from it. Since in general this is not the case, the contribution from the approximation error in the projected problem  can be estimated in terms of the eigenvalues $\lambda_j, \; j > k$. The projected equation (\ref{projected eq}) therefore can be seen as a non-linear extension of the linear spotlight inversion formula. By construction, this extension is different from the non-linear version in \cite{jaaskelainen2024projection}, which is based on projections onto subspaces orthogonal to the range of the Jacobian with respect to the nuisance parameter.

\subsection{Interpreting the clutter as Bayesian approximation error}\label{sec:clutter}

In this subsection we revisit spotlight inversion with the intent of interpreting the clutter term as a special case of approximation error. Recasting the linear inverse problem in the BAE setting, the matrix $\mA$ plays the role of the accurate model, the matrix $\mA_1$ that of the surrogate model, and the roles of $Z$ is played by $X_1$, so that 
\[
 f^*(x) = \mA x, \quad f(z) = \mA_1 z, \quad z = x_1.
\]
Following the Bayesian paradigm, we model all unknowns as random variable, and we define the approximation error as 
\[
 M = \mA X - \mA_1 Z = \mA_2 X_2,
\]
thus identifying the clutter term with the perturbation introduced by putting the spotlight on $X_1$. Next we consider how to estimate the probability density of $M$, first assuming that the prior of $X$ is Gaussian, then considering the non-Gaussian prior case.

Assume that a priori $X$ is modeled as a zero-mean Gaussian random variable, $X\sim{\mathcal N}(0,\mC_X)$, where $\mC_X$ is symmetric positive definite.  Partitioning the covariance matrix  as
\[
 \mC_X = \left[\begin{array}{cc} \mC_X^{11} & \mC_X^{12} \\ \mC_X^{21} & \mC_X^{22}
 \end{array}\right],
\]
where $\mC_X^{11} \in \R^{n \times n}$, and approximating $M$ by $\widetilde M$, which is independent of $X_1$ as in the EEM version of the BAE approach, we have
\begin{equation}\label{BME appr}
 \widetilde M \sim{\mathcal N}\left(0,\mA_2\mC_X^{22}\mA_2^\mT\right).
\end{equation}
Following the BAE paradigm, the approximate Maximum A Posteriori estimate $\widehat x_1$ is the minimizer of the quadratic expression
\[
 {\mathscr E}(x_1) = (b - \mA_1 x_1)^\mT(\mA_2\mC_X^{22}\mA_2^\mT + \mC_E)^{-1}(b-\mA_1 x_1) + x_1 \big(\mC_X^{11}\big)^{-1} x_1,
\]
or
\[
 \widehat x_1 = \big(\mA_1^\mT(\mA_2\mC_X^{22}\mA_2^\mT + \mC_E)^{-1}\mA_1^\mT + \big(\mC_X^{11}\big)^{-1}\big)^{-1}\mA_1^\mT
 (\mA_2\mC_X^{22}\mA_2^\mT + \mC_E)^{-1}b.
\]
The above formula underlines the attractiveness of the spotlight inversion approach in its implementational simplicity. Observe that under the assumption of the Gaussian prior, we may use the spotlight approach to the BAE approximation (\ref{BME appr}) and from the projector $\mP_k^\perp$ by using the eigenvectors of the covariance matrix of $\widetilde M$, it follows that
\[
 \mA_2\mC_X^{22} \mA_2^\mT = \sum_{j=1}^m \lambda_j u_j u_j^\mT.
\]
Observe that while the vectors $u_j$ form an orthonormal system in ${\mathscr R}(\mA_2)$, in general they do not coincide with the left singular vectors of the matrix $\mA_2$. Furthermore, the square roots $\sqrt{\lambda_j}$ of the eigenvalues  may decrease faster than the singular values of the matrix $\mA_2$, implying that they may provide a faster converging indicator of how large the cutoff index $k$ needs to be in order to have sufficient suppression of the clutter in the  projected model (\ref{projected eq}).

\subsection{Priorsketching}

In many application, it is not possible to find a Gaussian prior for $X$ that adequately  captures the underlying physics of the problem. A typical example where non-Gasussian priors are necessary is when the unknowns must be positive, as illustrated in the second computed examples in section ~\ref{sec:ComputedExamples}.

Without loss of generality, we may assume that $X$ has zero mean prior, implying that the mean of the approximation error $M = \mA_2 X_2$ vanishes.  The covariance matrix $\mC_M$ of the clutter $M$ is given by
\[
\mC_M = {\mathbb E} \big( (\mA_2 X_2) (\mA_2 X_2)^\mT \big) = \mA_2 {\mathbb E}(X_2 X_2^\mT) \mA_2^\mT. 
\]
If the covariance matrix of $X_2$ is not available, we estimate it numerically by drawing an independent sample of size $k$ from the prior for $X$, extract a corresponding sample for $Z = X_2$,  
\[
S= \left\{ x_2^{(1)}, \ldots, x_2^{(k)} \right\}.
\]
and form the sample matrix 
\[ \Omega = \frac1{\sqrt{k}}\left[ x_2^{(1)}, \ldots, x_2^{(k)} \right].
\]
The clutter covariance can now be estimated as
\begin{equation}\label{Cov M}
\mC_M \approx \mA_2 \Omega \Omega^\mT \mA_2^\mT.
\end{equation}
Observe that if the prior of $X_2$ is Gaussian as in the previous section, then the independent samples from the prior can be computed as
\[
 x^{(j)}_2 = \big(\mC_X^{22}\big)^{1/2} \xi^{(j)}, \quad \xi^{(j)}\sim {\mathcal N}(0,\mI_n),
\]
that is,
\begin{equation}\label{normal}
 \Omega = \big(\mC_X^{22}\big)^{1/2}\Xi, 
\end{equation}
where $\Xi\in\R^{n\times k}$ is a normal random matrix. 

To find few eigenvectors of the matrix $\mC_M$, there is no need to form the approximation (\ref{Cov M}), rather, we may compute the lean SVD  of the matrix $\mA_2\Omega \in\R^{m\times k}$,
\[
 \mA_2\Omega = \mU_k\mD_k\mV_k^\mT,
\]
and use the columns of $\mU_k$ to form the projection $\mP_k$ and its spotlight complement $\mP_k^\perp$. Alternatively, we may compute the QR-decomposition,
\[
 \mA_2\Omega = \mQ_k\mR_k,
\]
and use the columns of $\mQ_k$ to build the projectors. Either way, these may be seen as column sketching schemes of the matrix $\mA_2$. The difference from  the standard column sketching methods is that we draw the sketching matrix $\Omega$ from the prior of $X_2$, therefore we refer to it as {\em priosketching}. Observe that in the case when $X_2$ is Gaussian, the sketching coincides with normal column sketching (\ref{normal}) of the matrix $\mA_2\big(\mC_X^{22}\big)^{1/2}$. The interpretation of this approach is that rather than approximating the full range of the matrix $\mA_2$ by the orthonormal vectors, we concentrate on the part of the range in which the clutter is believed to live, in the light of the prior information. In particular, if the prior restricts $X_2$ to a lower dimensional subspace with dimension $k'<{\rm rank}(\mA_2)$, the range of $\mA_2\Omega$ is always at most $k'$-dimensional.



\subsection{Spotlighting as a limit case of BAE}

In spite of the many similarities, the solutions of the spotlight inversion and of the BAE compensated model reduction lead to  different solutions of the inverse problem. In this section we investigate on the differences between the two approaches and they affect the computed solutions. 

Under the additive white noise error assumption, in the Bayesian approximation error framework, we can write the surrogate inverse problem as 
\begin{equation}\label{centered}
 B -\mu = f(Z) + M + E,
\end{equation}
where 
\[
 M +  E \sim{\mathcal N}(0, \mC_M + \mI_m).
\]
Writing the approximation error covariance matrix  $\mC_M$ in terms of its eigenvalue decomposition, 
\[
 \mC_M + \mI_m  = \sum_{j=1}^m \lambda_j^2 u_ju_j^\mT + \mI_m
 =\sum_{j=1}^m (\lambda_j^2 + 1) u_ju_j^\mT. 
\]
and substituting it in the posterior model, we can write the data-dependent contribution in the Gibbs functional (\ref{MAP}) associated with likelihood model based on (\ref{centered}) as
\begin{equation}\label{Gibbs}
(b-\mu -f(z)^\mT(\mC_M + \mI_m)^{-1} (b-\mu -f(z)) = \sum_{j=1}^m \frac 1{\lambda_j^2 + 1}\left(u_j^\mT(b-\mu-f(x))\right)^2.
\end{equation}
If, for some $k<m$,  $\lambda_j^2 \ll 1$,
we can approximate (\ref{Gibbs}) as  
\begin{eqnarray*}
 &&(b-\mu -f(z)^\mT(\mC_M + \mI_m)^{-1} (b-\mu -f(z))\\
 &&\phantom{XXXXXX}\approx  \underbrace{\sum_{j=1}^k \frac 1{\lambda_j^2+1}\big(u_j^\mT(b-\mu-f(z))\big)^2}_{(I)}
 + \underbrace{\sum_{j=k+1}^m \big(u_j^\mT(b-\mu-f(z))\big)^2}_{(II)}
\end{eqnarray*}
The spotlight inversion scheme projects (\ref{centered}) onto the orthocomplement of the $k$-dimensional approximation of the modeling error, completely neglecting the term (I), accounting for the fidelity term only on the orthogonal complement of it. The BAE model, on the other hand, which accounts for both (I) and (II), weights less heavily data mismatch in the subspace of the approximation error, effectively assigning to the projection of the discrepancy in each eigendirection of the approximation error covariance matrix a weight that is inversely proportional to the corresponding eigenvalue. 
Effectively, the spotlighting approach can be seen as a limit case of the BAE approach, where the $k$ dominant eigenvalues are sent to infinity.  

\section{Computed examples}\label{sec:ComputedExamples}

We elucidate the discussion of the previous sections with two computed examples. The first one is a linear inverse problem arising in X-ray tomography, and the second one is the non-linear electrical impedance tomography (EIT) problem with geometric model uncertainties.

\subsection{X-ray tomography}

In the first computed example, we consider fanbeam X-ray tomography in two dimensions. In the experiment, the 
domain modeled as a unit square, $\Omega = [0,1]\times [0,1]$, is illuminated by a fanbeam X-ray source from several directions, and the attenuation of the radiation traversing the domain is recorded on the opposite side. Based on the Beer-Lambert law, the data for each illumination angle can be modeled as integrals of the attenuation coefficient along the beam lines through the domain.
The unknown of interest is the distributed
attenuation parameter in $\Omega$.

\begin{figure}
 \centerline{\includegraphics[width=12cm]{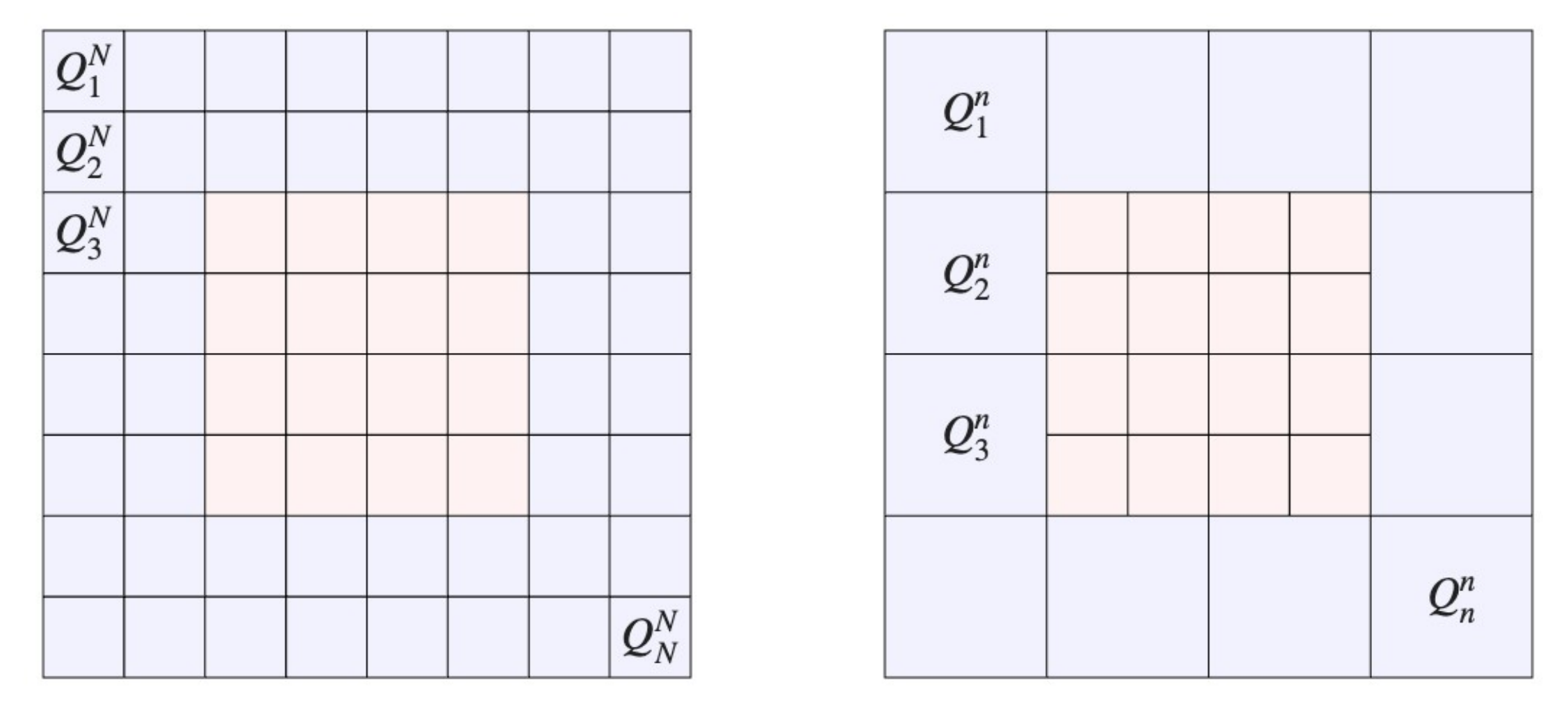}}
 \caption{A schematic picture of the coarsening. Here, the pink square is the region of interest $D$, where we keep the fine discretization.}\label{fig:pixels}
\end{figure}

To define a discrete model,
we represent the domain as a union of $N$ fine pixels obtained by dividing the unit interval into $\sqrt{N}\in\N$ subintervals, see Figure~\ref{fig:pixels} and denote the $j$th pixel by $Q^N_j$, $1\leq j\leq N$,.

Suppose that we are mostly interested in a detail of the image, referred to as the {\em spotlight domain} $D\subset\Omega$. For the sake of simplicity, we select $D$ to be a rectangular domain consisting of $d$ pixels. 
We create a coarse model of lower computational complexity by aggregating the pixels outside the spotlight domain into bigger pixels. Let $n$ denote the number of pixels in the coarse model, which is the union of the fine pixels inside $D$ and the lumped pixels outside $D$. We denote the pixels in the coarse model by $Q^n_j$, $1\leq j\leq n$.

Consider a fine grain image $x^N\in\R^N$. To map this image to a coarse grain image by averaging the pixel values of the fine image over the coarse pixels, we define a matrix $\mP\in\N^{n\times N}$
such that
\begin{eqnarray*}
 \mP_{kj} & = & 1\quad \mbox{ if and only if $Q^N_j\in Q^n_k$,} \\
 & = &  0 \quad \mbox{ otherwise}.
\end{eqnarray*}
 Multiplying $x^N$ by $\mP$ adds up the values of the pixels of the finer model inside each pixel of the coarse model, thus leaving the pixel values in the domain $D$ unchanged. To ensure that both fine and coarse images have the same intensity, we define the scaling matrix  
 \begin{equation}\label{reduce image}
  \mW = {\rm diag}\left(\sum_{j=1}^N \mP_{1,j}, \ldots, \sum_{j=1}^N \mP_{n,j}\right)\in\R^{n\times n}.
 \end{equation}
 and let 
 \[
  x^n = \mW^{-1} \mP x^N,
  \]
  so that the value of each coarse pixel is the average of the values in the fine pixels defining it.

We define the tomography matrix $\mA^N \in\R^{m\times N}$, where $m$ is the number of illumination rays in the model, by assigning the value of the element in the $\ell$th row and $j$th column to be the length of the intersection of the $\ell$th ray and the $j$th pixel, i.e., 
\[
 \mA^N_{\ell,j} = |L_\ell\cap Q_j^N|, \quad \mbox{$L_\ell = $ the $\ell$th ray passing through the object.}
\] 
Similarly, we define the tomography matrix  $\mA^n \in\R^{m\times n}$ for the coarse model in the same manner by setting
\[
 \mA^n_{\ell,j} = |L_\ell\cap Q_j^n|.
\] 
Observe that since each coarse pixel is the aggregation of fine pixels, 
\[
 \mA^n_{\ell,j} = |L_\ell\cap Q_j^n| = \sum_{Q_k^N\subset Q_j^n}   |L_\ell\cap Q_k^N|  = \sum_{k=1}^N \mP_{j k}  \underbrace{ |L_\ell\cap Q_k^N| }_{=\mA^N_{\ell k}},
 \] 
 which can be expressed in matrix form as 
 \begin{equation}\label{reduce matrix}
  \mA^n = \mA^N \mP^\mT, 
 \end{equation}
thus establishing the relation between the tomography matrices of the fine and coarse model. 

We illustrate the performance of the two model reduction algorithms for computed tomography using the lotus root data set described and downloadable in \cite{bubbaOA}. The fine scale pixel map is of size $N = 128\times 128 = 16\,384$, and the data consist of a sinogram of $120$ projection angles over the interval $[0,2\pi]$. Since each projection comprises $429$ rays passing through $\Omega$, the dimension of the data space is $m = 120\times429 = 51\,480$. To set the noise level, we estimated standard deviation of the presumably Gaussian additive scaled white noise, obtained by considering the data corresponding to rays that miss the target.

\begin{figure}[ht!]
\centerline{
\includegraphics[width=16cm]{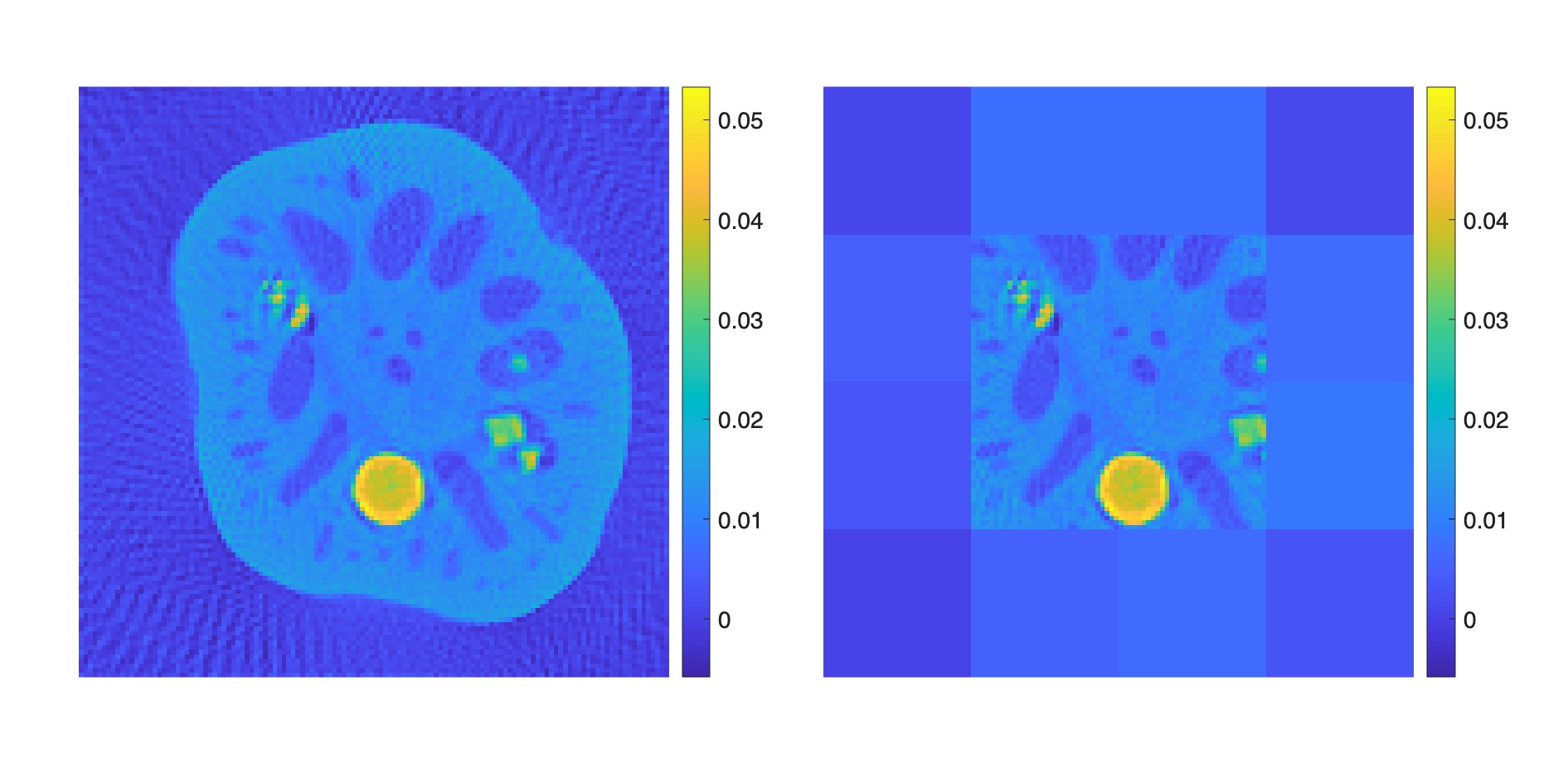}
}
\caption{On the left, we show the least squares solution of the problem $\mA^N x^N = b$ solved by using the Krylov subspace LSQR iterative algorithm, regularized by stopping the iterations when the discrepancy reaches the estimated noise level. The image consists of $128\times 128$ pixels. The image shown on the right was obtained by subdividing the image on the left into $4\times 4$ subdomains, each comprising $32\times 32$ pixels, assigning to the pixel in the border subdomains the average value of the finer are averaged, while pixel values in the four center subdomains are intact. We consider this image as the reference reconstruction and denote it by $x^n_{\rm ref}$.}\label{fig:full data}
\end{figure}

We start by computing an estimate of the image $x^N$ by solving the fine grain problem
\begin{equation}\label{fine grain}
\mA^N x^N = b
\end{equation}
in the least squares sense, using the Krylov subspace based LSQR algorithm, stopping the iterations when the norm of the discrepancy reaches the noise level.
Observe that for this mildly ill-posed inverse problem, the stopping criterion provides a regularization that follows the Morozov discrepancy principle, see, e.g.,  \cite{kilmer2001choosing}.
The reconstruction is shown in the left panel of Figure~\ref{fig:full data}. To create a coarse scale model, we subdivide first the image $x^N$ into $4\times 4$ large subdomains, each one comprising $32\times 32$ pixels, and aggregate the   pixels in each large subdomain, except for the four  $32\times 32$ subdomains at the center of the image, which constitute the region of interest where we want to solve the problem with high accuracy. 
This region of interest, referred to as the spotlight domain $D$, therefore consists of $d = 64\times 64 = 4\,096$ fine pixels. 

The dimension of the coarse model, which consists of the 12 coarse pixels on the border and the finely discretized inner region of interest is $n = 4\,096 + 12 = 4\, 108$. The right panel in Figure~\ref{fig:full data} shows the lumped solution $x^n = \mS^{-1}\mP x^N$.
Since the data where obtained experimentally, thus no ground truth is provided, we use the image $x^n$ built as described above as our reference solution for testing the performance of the different methods, and we denote it by   $x^n_{\rm ref}.$

\begin{figure}[ht!]
\centerline{
\includegraphics[width=16cm]{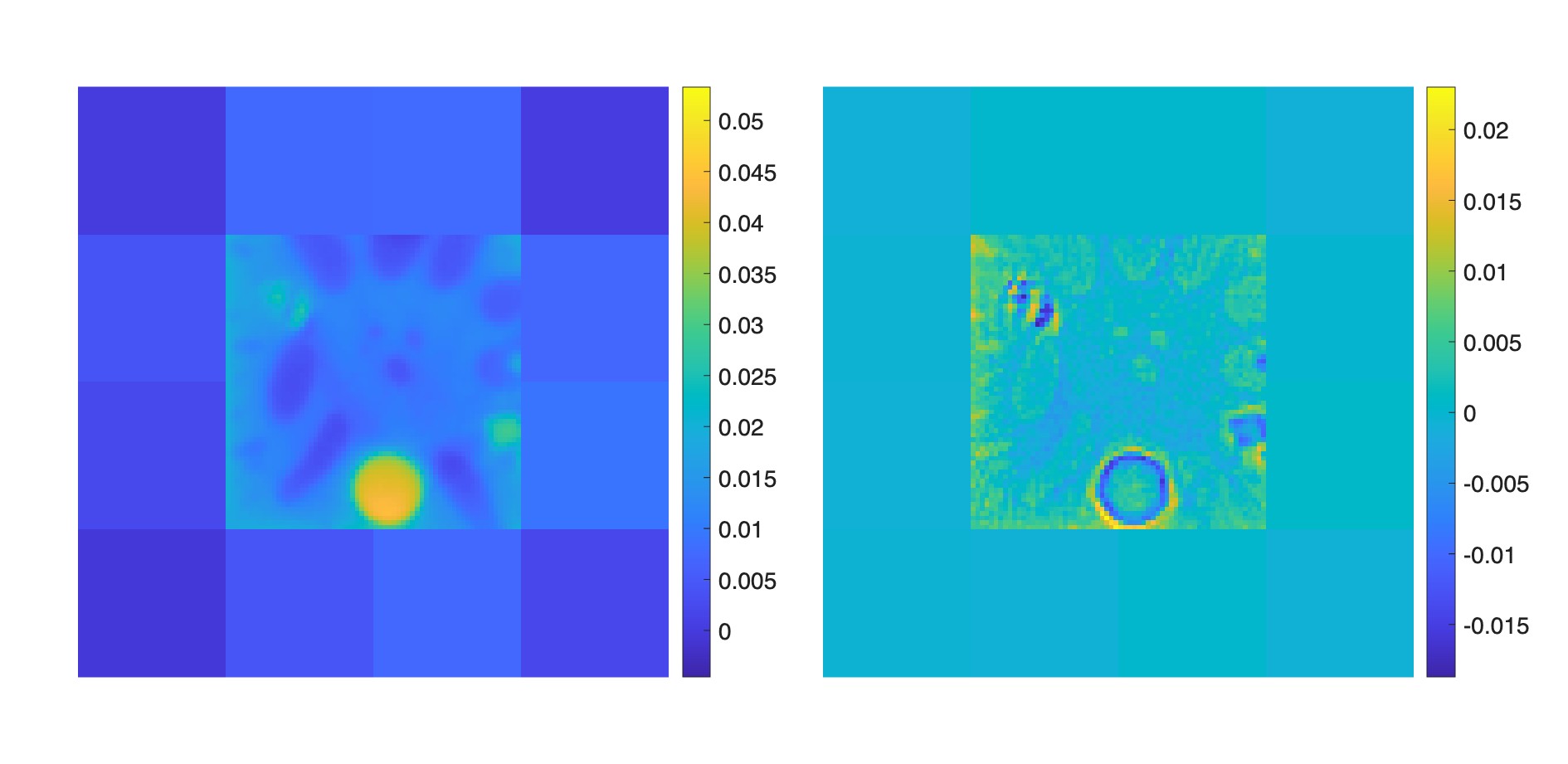}
}
\caption{On the left we show the reconstruction $x^n_{\rm naive}$ obtained by ignoring the approximation error and solving the equation $\mA^n x^n = b$ by the LSQR algorithm, stopping the iterations when the discrepancy reaches the noise level. On the right, we show the difference $x_{\rm naive} - x_{\rm ref}$.}\label{fig:no mod error}
\end{figure}

To demonstrate the importance of accounting for the error introduced by a model reduction, we look at what happens when the approximation error is ignored. To that end, we estimate $x^n$ by solving in the least squares sense the equation
\begin{equation}\label{coarse grain}
\mA^n x^n = b,
\end{equation}
again using the Krylov subspace based LSQR solver, regularizing the problem by stopping the iterations when the discrepancy reaches the noise level. The reconstruction obtained in this manner, and its deviation from the reference solution $x^n_{\rm ref}$, are shown in Figure~\ref{fig:no mod error}. We refer to this solution as the na\"{\i}ve solution, and denote it by $x^n_{\rm naive}$. Neglecting the approximation error causes a visible blurring of the spotlight domain and introduces a boundary halo along the edges of the region of interest. The deviation from the reference image is significant in amplitude, locally arriving to values of approximately one third of the dynamic range of the reference image.

\begin{figure}[ht!]
    \centerline{
    \includegraphics[width=16cm]{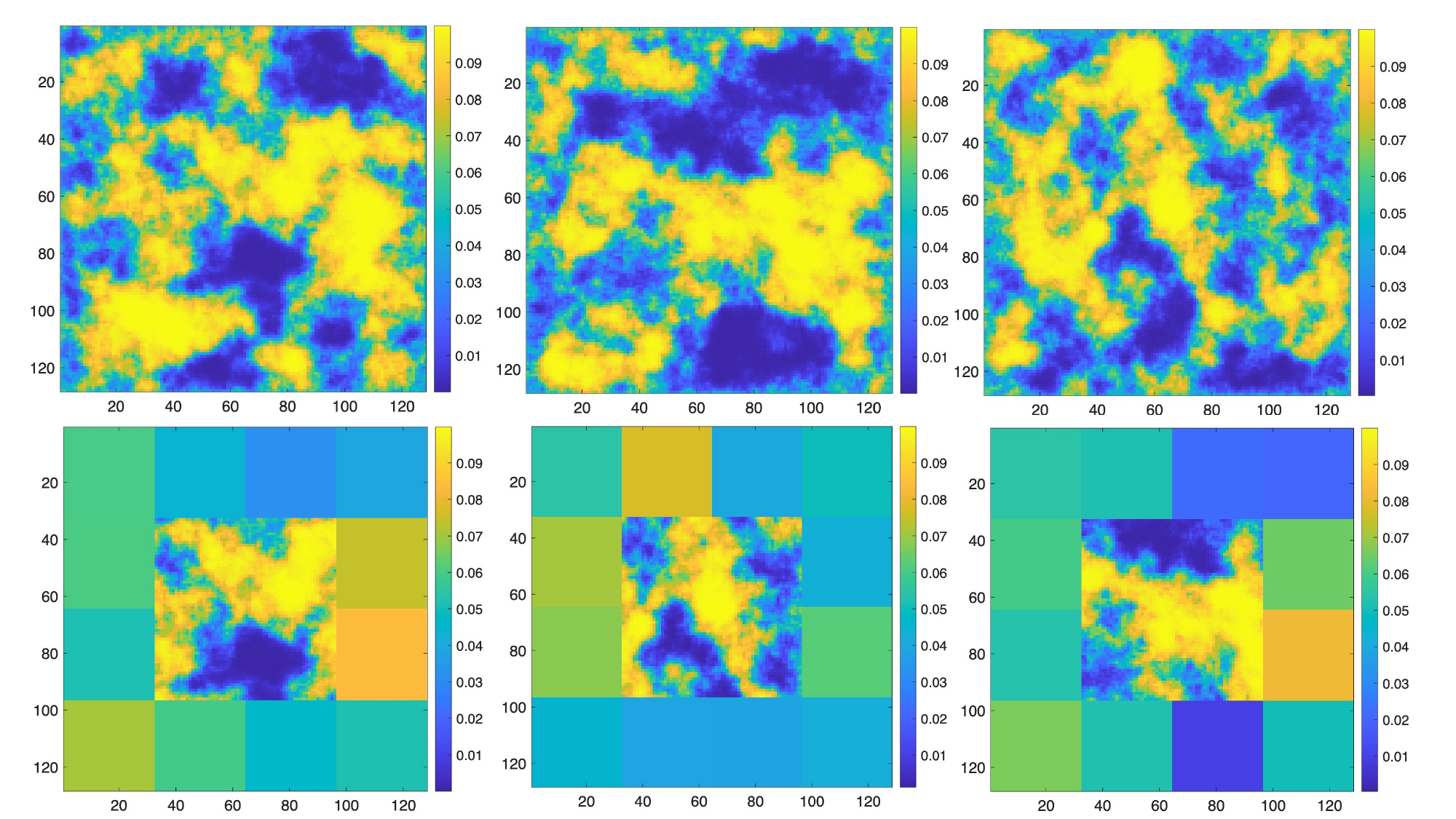}
    }
    \caption{Three draws from the prior model for the fine grain variable (top row) and the corresponding coarse counterparts.}\label{fig:tomo draws}
\end{figure}

Having established the importance of including the approximation error, we consider next the Bayesian error model. To this end, a prior model for generating the approximation error sample needs to be defined.
To ensure that the absorption coefficient is non-negative, we define a non-Gaussian prior model, writing
\begin{equation}\label{random1}
 X = \gamma\, \Sigma_{\xi_0,\alpha}\big(\Xi\big),
\end{equation}
where $\Xi$ is a Gaussian random variable modeled as
\begin{equation}\label{radom2}
(-\Delta + \lambda^{-2} \mI)\,\Xi = W\sim{\mathcal N}(0,\mI), \quad \lambda = \mbox{correlation length},
\end{equation}
and $\Sigma_{\xi_0,\alpha}:\R^N\to \R^N$ is component-wise sigmoid function,
\begin{equation}\label{random3}
\big(\Sigma_{\xi_0,\alpha}(\xi)\big)_j = \frac{1}{1+e^{\alpha(\xi_0 -\xi_j)}},
\end{equation}
The parameter $\xi_0$ marks where the sigmoid assumes its mid value $1/2$, while $\alpha$, proportional to the slope at the half value, controls the contrast in the image, as the parameter determines how sharply the function changes from zero to one as $\xi$ passes the half value. The parameter $\gamma$ is a scaling factor that makes it possible to tune the dynamical range of the prior.
Figure~\ref{fig:tomo draws} shows three random draws for the prior, as well as the corresponding projections to the coarse grid.

Once a sample of the modeling error has been computed, it can be used to setup the reduced model inversions, the BAE compensation and the spotlight projection model.

\subsubsection{Computing the solution of the BME reduced model}

Consider the approximate problem
\[
 b = \mA^n x^n + M + \varepsilon, 
\]
where $M$ is the approximation error and that $\varepsilon$ is scaled white noise with variance $\sigma^2$. Given a sample of the approximation error
\[
 \{m^{(1)}, m^{(2)}, \ldots,m^{(k)}\},
\]
we use it to compute the empirical mean and covariance 
\[
 \mu = \frac 1k\sum_{j=1}^k m^{(j)}, \quad 
 \mC = \frac 1k \sum_{j=1}^k \mu^{(j)}(\mu^{(j)})^\mT, \quad \mu^{(j)} = m^{(j)}-\mu.
\]
To highlight the fact that Bayesian approximation error compensation and spotlight inversion do not require the inverse problem to be solved in the Bayesian framework, the solution of the reduced model is computed by minimizing the log-likelihood using an iterative linear solver regularized by stopping the iterations as soon as the norm discrepancy is at the level of the norm of the observation error. 
In the case of the BAE reduced model, the objective function to minimize is
\begin{equation}\label{obj}
 {\mathscr E}(x^n) = (b - \mu - \mA^n x^n)^\mT (\mC + \sigma^2\mI)^{-1}(b - \mu- \mA^n x^n).
\end{equation}
If the dimensionality of the data space is small, it is possible to use direct methods relying on the factorization of the precision matrix. When that is not feasible, or the precision matrix is not available, it is possible to arrange the computation required by the iterative solvers taking advantage of the fact that $\mC$ is a low-rank approximation of the approximation error covariance. Writing
\begin{equation}\label{S matrix}
 \mC = \mS \mS^\mT, \quad \mS = \frac1{\sqrt{k}}\left[\begin{array}{ccc} \mu^{(1)} & \cdots &\mu^{(k)}\end{array}\right],
\end{equation}
it follows from the Sherman-Morrison-Woodbury identity that
\begin{equation}\label{SMW} 
 (\mC + \sigma^2 \mI)^{-1} = \frac 1{\sigma^2}\left( \mI - \mS (\sigma^2 \mI_k + \mS^\mT\mS)^{-1}\mS^\mT\right) 
 = \frac 1{\sigma^2}\left( \mI - \mK \mK^\mT\right),
\end{equation}
with $\mK =\mS\mR^\mT \in\R^{m\times k}$ where $\mR$ is the Cholesky factor of a $k\times k$ matrix, i.e., 
\[
 \mR^\mT\mR = (\sigma^2 \mI_k + \mS^\mT\mS)^{-1}.
\]
Rewriting  (\ref{obj}) in view of (\ref{SMW}), the expression to be minimized becomes 
\[
  {\mathscr E}(x^n) = \frac1{\sigma^2}(b - \mu - \mA^n x^n)^\mT ((\mI - \mK\mK^\mT)(b - \mu- \mA^n x^n),
\]
and the minimizer of this quadratic form satisfies
\[
 (\mA^n)^\mT(\mI - \mK\mK^\mT)\mA^n x =(\mA^n)^\mT(\mI - \mK\mK^\mT)(b-\mu).
\]
When $m$ is large and the matrix $\mK\mK^\mT$ cannot be formed explicitly, the evaluation of the left hand side can be implemented in a matrix-free manner provided that we have access to a routine performing matrix-vector products with $\mK$ and $\mK^\mT$, 
\[
 x\mapsto \mA^n x\mapsto \mK^\mT\mA^n x\mapsto
 \mK\mK^\mT\mA^n x \mapsto (\mA^n)^\mT\mK\mK^\mT\mA^n x.
\]
In our computed example, the solution of the BAE reduced problem was computed using the Conjugate Gradient algorithm in a matrix-free manner, regularizing the problem by stopping the iterations as soon as the discrepancy reaches the noise level.

\begin{figure}[ht!]
\centerline{
\includegraphics[width=16cm]{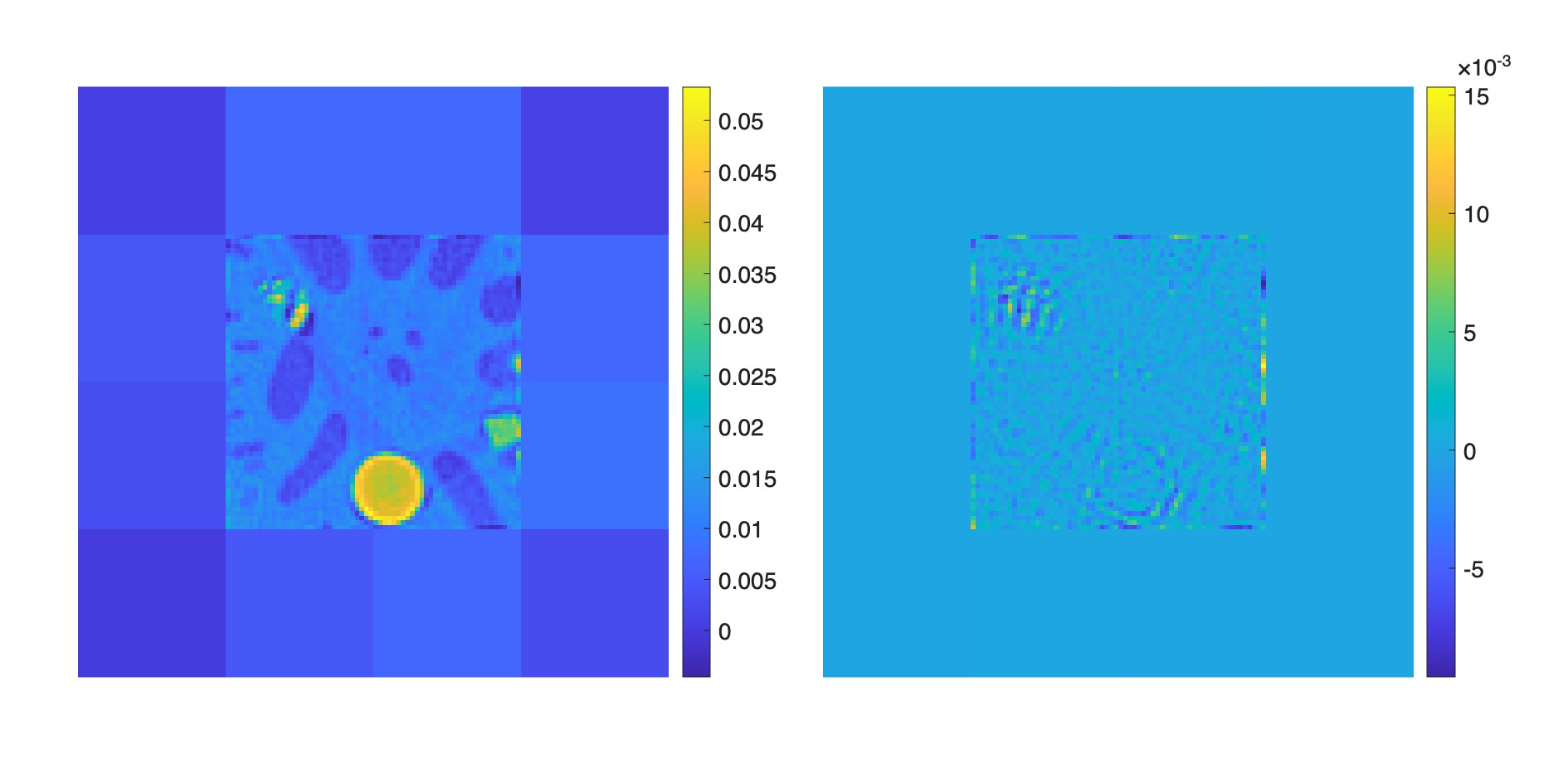}
}
\caption{On the left, we show the reconstruction $x^n_{\rm BME}$ obtained by using the BAE reduced model. The solution is calculated by using the Conjugate Gradient solver regularized by terminating the iterations as soon as the discrepancy norm reached the noise level.
On the right, we show the difference $x_{\rm BAE} - x_{\rm ref}$.}\label{fig:enhanced error}
\end{figure}
To test this approach, we generate a sample of 250 draws from the prior described in the previous section. The result, together with the discrepancy from the reference reconstruction is shown in Figure~\ref{fig:enhanced error}. The reconstruction is superior to that obtained ignoring the approximation error, as expected, the discrepancy from the reference image dropping about an order of magnitude. Numerical tests show that reducing the number of draws decreases slightly the quality, in particular, causing a slight boundary effect around the spotlight region.

\subsubsection{Solution by the spotlight projection}

Consider next the spotlight inversion, where we represent the approximation error by its eigenvector expansion (\ref{eig expansion}) and approximate the orthogonal projector by the priorsketching of the approximation error. We have the projected equation
\[
 \mP^\perp_k(b -\mu) = \mP^\perp_k \mA^n z^n 
 +\sum_{j=k+1}^m\sqrt{\lambda_j}\xi_ju_j  +\mP^\perp_k \varepsilon,
\]
and assuming that the projected approximation error is insignificant, we solve the projected least squares problem
\[
 \mbox{minimize $\|\mP^\perp_k(b -\mu) - \mP^\perp_k \mA^n z^n \|^2$}
\]
by using the Krylov subspace based LSQR algorithm for linear least squares with stopping by the Morozov discrepancy principle. We remark that while the approximation error sample was generated by using a prior model, the solution strategy of the projected problem need not align with the Bayesian paradigm, the argument being purely linear algebraic. 

\begin{figure}[ht!]
\centerline{
\includegraphics[width=16cm]{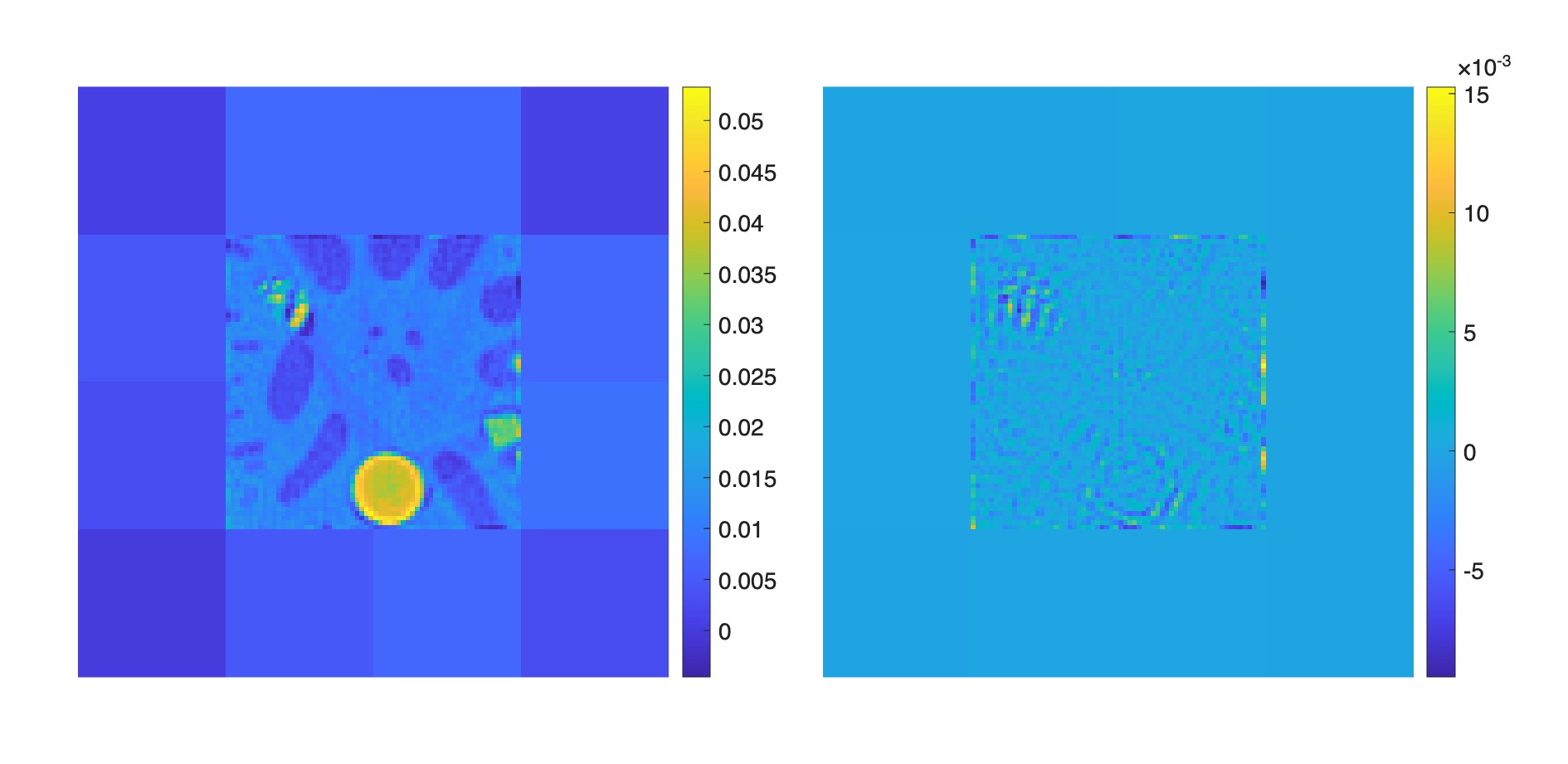}
}
\caption{The spotlight solution $x^n_{\rm spot}$, the range of the approximation error being estimated by using the prior sketching, and the difference $x^n_{\rm spot} - x^n_{\rm ref}$.}\label{fig:spotlight}
\end{figure}

Using the same sample of 250 draws of the approximation error, we form the matrix $\mS$ as in (\ref{S matrix}), and compute the vectors $u_j$ as the right singular vectors of the matrix $\mS$. The result of the projected least squares solution, shown in Figure~\ref{fig:spotlight}, is remarkably similar to the solution of the BAE reduced model, however the spotlight implementation does not require at all the approximation error covariance matrix or its inverse. 





\subsection{Electrical impedance tomography}

In the second example, we consider the inverse problem of electrical impedance tomography (EIT). Given a bounded domain $\Omega$ with connected complement, the goal is to estimate the electric conductivity $\sigma(x)>0$ inside $\Omega$ by current-voltage measurements on the boundary $\partial\Omega$. We assume that $L$ contact electrodes $e_\ell \subset\partial\Omega$ are attached to the boundary, $L$ currents $I_\ell$ are injected through them into or out of the body, and the voltages corresponding $V_\ell$ are measured on each electrode. Since the currents need to satisfy Kirchhoff's law of conservation of charge,
\[
 \sum_{\ell=1}^L I_\ell = 0,
\]
there are $L-1$ linearly independent current patterns, that can be chosen to form a complete current frame $\{I^1,\ldots,I^{L-1}\}$.
We denote the corresponding voltage patterns by 
$\{V^1,\ldots,V^{L-1}\}$, and assume that the ground voltage is chosen so that the condition
\[
 \sum_{\ell=1}^L V^k_\ell = 0
\]
holds for every $V^k\in\R^L$. We set up the forward model using the complete electrode model \cite{somersalo1992existence}, assuming that the voltage potential $u$ satisfies the equation
\[
 \nabla\cdot\big(\sigma\nabla u \big) = 0 \mbox{ in $\Omega$,}
\]
with the boundary conditions
\[
 \int_{e_\ell} \sigma\frac{\partial u}{\partial n} dS = I_\ell, \quad \sigma\frac{\partial u}{\partial n}\bigg|_{\partial\Omega\setminus\cup e_\ell} = 0,
\]
and
\[
 \left(u + z_\ell\sigma \frac{\partial u}{\partial n}\right)\bigg|_{e_\ell} = V_\ell,
\]
where $z_\ell>0$ is the contact impedance of the $\ell$th electrode. The ensuing boundary value problem is well-posed, and has a unique solution
\[
 (u,V) \in H^1(\Omega)\times \R^L_0, \quad \R^L_0 = \big\{V\in\R^L\mid \sum_{\ell =1}^L V_\ell = 0\big\},
\]
that can be computed using a finite element implementation; see, e.g.,  \cite{calvetti2025sparsity}.

It is well-known (see, e.g., \cite{kolehmainen2013recovering,kolehmainen2008electrical,hyvonen2017compensation,candiani2019computational}) that the EIT inverse problem is sensitive to geometric uncertainties about the shape of the boundary of the domain and the position of the electrodes. The BAE approach makes it possible to mitigate the reconstruction artifacts arising from inaccurate modeling of the boundary shape. In this computed example, we show that the projection method provides an efficient and computationally simple implementation of the approximation error approach. For comparison, we refer in particular to the articles \cite{nissinen2011reconstruction,kaipio2013approximate}, where the problem was addressed using the BAE approach with marginalization with respect to the geometric uncertainty, and low rank approximations of the error covariance matrix were used in computations.

Consider the EIT experiment in which a fixed set of current patterns are applied on the electrodes, and the corresponding voltages are measured on the surface of the body $\Omega$ whose exact shape is poorly known. To solve the inverse problem, we define a reference domain $\Omega_0$, and denote by $\Omega_\omega$
random perturbations of the reference shape, where $\omega$ is a sample from a probability space. We assume that for each $\omega$, there is a diffeomorphism $T_\omega:\Omega_0\to\Omega_\omega$, and the conductivities $\sigma$ in $\Omega_0$ and  $\sigma_\omega$ in $\Omega_\omega$ are related through
$\sigma_\omega(x) = \sigma(T^{-1}_\omega(x)) = T_{\omega,\#}\sigma(x)$, where $T_{\omega,\#}$ is the push-forward map. Thus, if the conductivity is estimated using the reference domain $\Omega_0$, it can be mapped to the domain $\Omega_\omega$ in an unambiguous way. 

\begin{figure}[ht!]
    \centerline{
    \includegraphics[width=16cm]{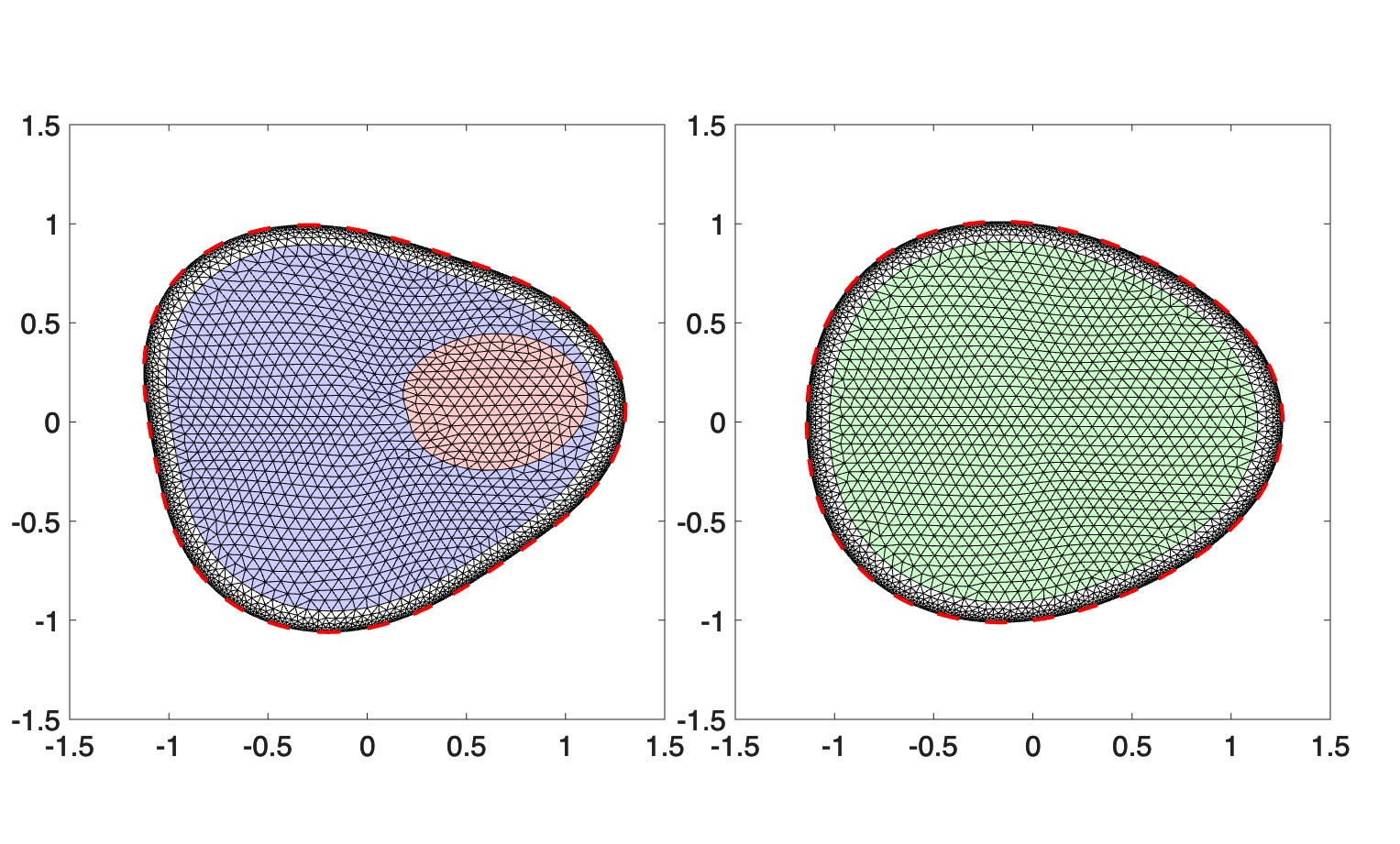}
    }
    \caption{The left panel shows the generative model. The conductivity is set to the background value $\sigma = \sigma_0$ outside the red inclusion, where the conductivity is $\sigma_1 = 3$. On the right, the standard shape is shown with the reconstruction mesh. Outside the green region, the conductivity is assumed to be known and fixed at the value $\sigma_0$. The 32 electrodes are marked by the red line segments along the boundaries.}\label{fig:meshes}
\end{figure}

For the sake of definiteness, let the set of applied current patterns constitute a frame of $L-1$ linearly independent patterns, and after stacking the the $L-1$ voltages in a vector $V$ of length $m = L(L-1)$, we write the model as
\begin{equation}\label{EITmodel}
 V = F_\omega(\sigma_\omega) + e,
\end{equation}
where $F_\omega$ is a randomized forward map corresponding to the domain $\Omega_\omega$. To estimate the conductivity using the reference model
$F_0$ based on the reference domain $\Omega_0$, we write
\begin{equation}\label{forward}
 V = F_0(\sigma) + \big\{F_\omega(\sigma_\omega) -
F_0(\sigma)\big\} + e = F_0(\sigma) + M_\omega(\sigma) + e,
\end{equation}
where the approximation error term is
\begin{equation}\label{modeling error}
 M_\omega(\sigma) =
 F_\omega(\sigma_\omega) -
F_0(\sigma) = F_\omega\big(T_{\omega,\#}\sigma\big) - F_0(\sigma).
\end{equation}
Observe that since in the current application we do not know the true shape of the boundary, the accurate forward map is also unknown. In fact, if the shape were known, it would be used in the inversion rather than resorting  to the reference shape. Nonetheless, we assume implicitly that the prior distribution used to estimate the approximation error can represent all the possible domain shapes, including the true, or generative shape. In Figure~\ref{fig:meshes}, the generative shape with an inclusion and the reference shape with the reconstruction mesh are shown.

The reference shape shown in the panel on the right is obtained by mapping a circular unit disk
to the domain $\Omega_0$ through the radial coordinates $(r,\theta)$ 
of the unit disc as
\[
 (x,y) \mapsto (r,\theta) \mapsto (\rho(r,\theta),\theta) \mapsto (x',y') =(1.1\,\rho(r,\theta)\cos\theta,\rho(r,\theta)\sin\theta), 
\]
where 
\begin{equation}\label{standard shape}
\rho(r,\theta) = r(1 + 0.05\,\cos 3\theta).
\end{equation}
The generative model, on the other hand, is a realization of a randomized radial transformation,
\begin{equation}\label{random shape}
\rho(r,\theta,\omega) = r\big(1 + 0.1 \xi(\omega)\,\cos 3\theta + 
0.1(\nu(\omega)-1/2)\sin 3\theta\big), \quad\xi,\nu\sim{\rm Uniform}([0,1])
\end{equation}
Since both the reference domain and the domains $\Omega_\omega$ are diffeomorphic with the unit disc, the transformations generate the diffeomorphism $T_\omega:\Omega_0\to \Omega_\omega$, while its explicit form is not of interest here.

\begin{figure}[ht!]
    \centerline{
    \includegraphics[width=18cm]{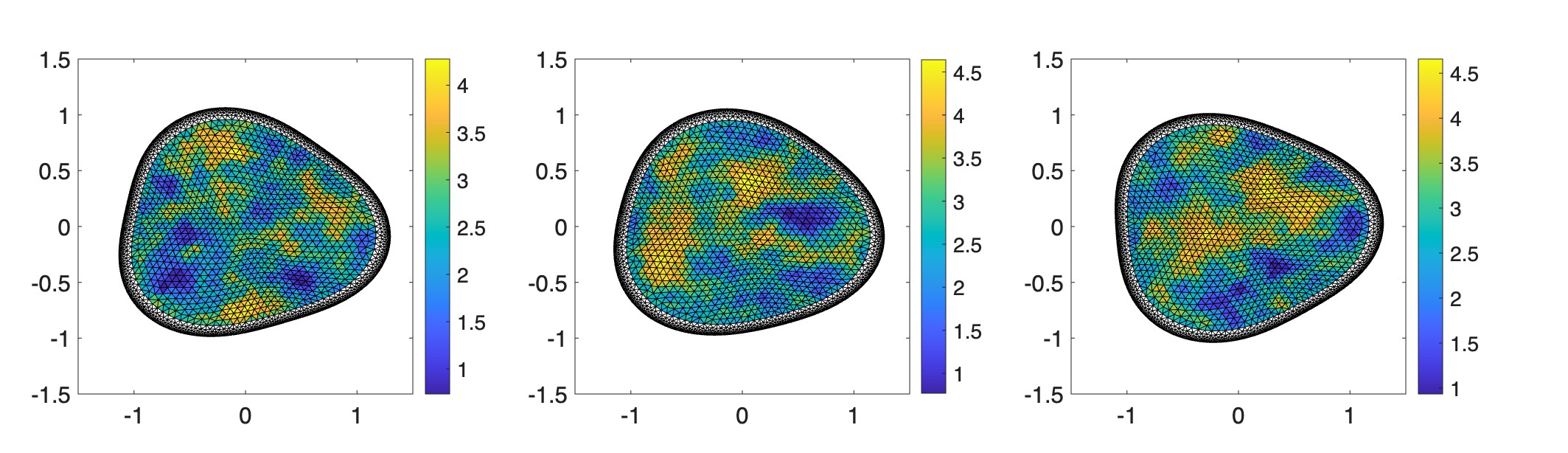}
    }
    \caption{Three random conductivity structures inside the corresponding randomly varying domains.}\label{fig:conductivity draws}
\end{figure}

To generate a sample of the approximation errors, we draw a sample of domains $\Omega_{\omega_j}$, $1\leq j\leq k$, from the model (\ref{random shape}), and generate corresponding pairs of random conductivities $(\sigma^{(j)},T_{\omega^j,\#}\sigma^{(j)})$, $1\leq j\leq k$, by generating the conductivities in the unit disc and  mapping them by the push-forward maps corresponding to (\ref{standard shape}) and (\ref{random shape}), respectively. The computation of the approximation error sample $M_{\omega_j}(\sigma^{(j)})$ from (\ref{modeling error}) requires solving the corresponding two FEM problems. We draw the conductivities using a prior model analogous to the one used for the tomography example, assuming that the conductivity is piecewise constant in the element basis over the unit disc. More precisely, we assume that in the unit disc, $\sigma = 1$ outside the disc of radius $r=0.9$, which corresponds to the green region in the standard domain of Figure~\ref{fig:meshes}, while inside the disc, the conductivities are drawn from the model (\ref{random1})-(\ref{random3}), with the Laplacian  approximated by the discrete Laplacian in the element basis, $\lambda = 5$, $\alpha = 3$, and $\gamma = 5$. Figure~\ref{fig:conductivity draws} shows three conductivity structures plotted in the corresponding randomly deformed domains.

Once a sample of size $k$ of the approximation error vectors has been generated,
\[
 M^j = M_{\omega_j}(\sigma^j) =
 F_{\omega_j}\big(T_{\omega_j,\#}\sigma^j\big)- F_0\big(\sigma^j\big), \quad 1\leq j\leq k,
\]
we write the approximation
\[
 M \approx \mu + \sum_{j=1}^k \sqrt{\lambda_j} \xi_j u_j,
\]
where
\[
 \mu = \frac 1k\sum_{j=1}^k M^j.
\]
The vectors $u_j$, which are the left singular vectors of the matrix 
\[
 \mM_c = \frac 1{\sqrt{k}}\left[\begin{array}{ccc} M^1-\mu & \cdots &M^k-\mu\end{array}\right],
\]
define the orthogonal projector pair 
\[
 \mP = \sum_{j=1}^k u_ju_j^\mT, \quad \mP^\perp = \mI_m - \mP
\]
onto the subspace ${\rm span}\{u_1,\ldots,u_k\}$ and its orthogonal complement, respectively. The projected problem that we consider is 
\begin{eqnarray*}
 \mP^\perp(V - \mu) &=& \mP^\perp F_0(\sigma) + \mP^\perp\big( M_\omega(\sigma) -\mu\big) + \mP^\perp e  \\
 &\approx& \mP^\perp F_0(\sigma)  + \mP^\perp e,
\end{eqnarray*}
assuming that after the projection by $\mP^\perp$ the contribution of the approximation error is at the level or below that of the observation noise. The conductivity is then estimated by solving a sequence of regularized least squares problems. More precisely, to guarantee that the conductivity is positive, we parameterize it as 
\[
 \sigma_j = \sigma_0 e^{x_j},
\] 
and, for the sake of simplicity, we assume that $\sigma_0=1$. 
We introduce the notation
\[
 G(x) = F_0(e^x),
\]
and write a Tikhonov-regularized objective function,
\[
 {\mathcal E}(x) = \|\mP^\perp G(x) -\mP^\perp(V-\mu)\|^2  + \delta \|\mL x\|^2,
\]
where the functional $\mL$ is a discrete approximation of 
the operator $-\Delta + \lambda^{-2}$, i.e., a Whittle-Mat\'{e}rn type smoothness penalty, $\lambda$ being a correlation length parameter. We remark that the regularization operator is not related to the prior used to generate the modeling error sample, other than for the fact that we use the same correlation length parameter. The regularization operator here is chosen purely for computational convenience. By linearizing $G$  around the current value $x_c$ of $x$, we approximate the first term of the objective function as
\[
    \|\mP^\perp G(x_c+\delta x) - \mP^\perp (V-\mu) \|^2 \approx \|
    \underbrace{\mP^\perp DG(x_c)}_{\mA}\delta x - \underbrace{\mP^\perp( V-\mu) - \mP^\perp G(x_c)}_{=r} \|^2 \\
    = \| \mA\delta x - r\|^2
\]
so that the local minimization problem reduces to solving the least squares problem   
\begin{equation}\label{regularized}
\|\mA\delta x - r\|^2 + \| \mL(x_c + \delta x)\|^2 = \left\|\left[\begin{array}{c} \mA \\ \mL\end{array}\right]\delta x - \left[\begin{array}{c} r \\ -\mL x_c\end{array}\right]
\right\|^2,
\end{equation}
whose solution solves the associated normal equations
\[
 \left(\mA^\mT\mA + \mL^\mT\mL\right)\delta x = \mA^\mT r - \mL^\mT\mL x_c.
\]
In our numerical simulations, we perform three steps of the above Gauss-Newton iterations. After the three steps, the solution does not change significantly.

\begin{figure}[ht!]
\centerline{
\includegraphics[width=16cm]{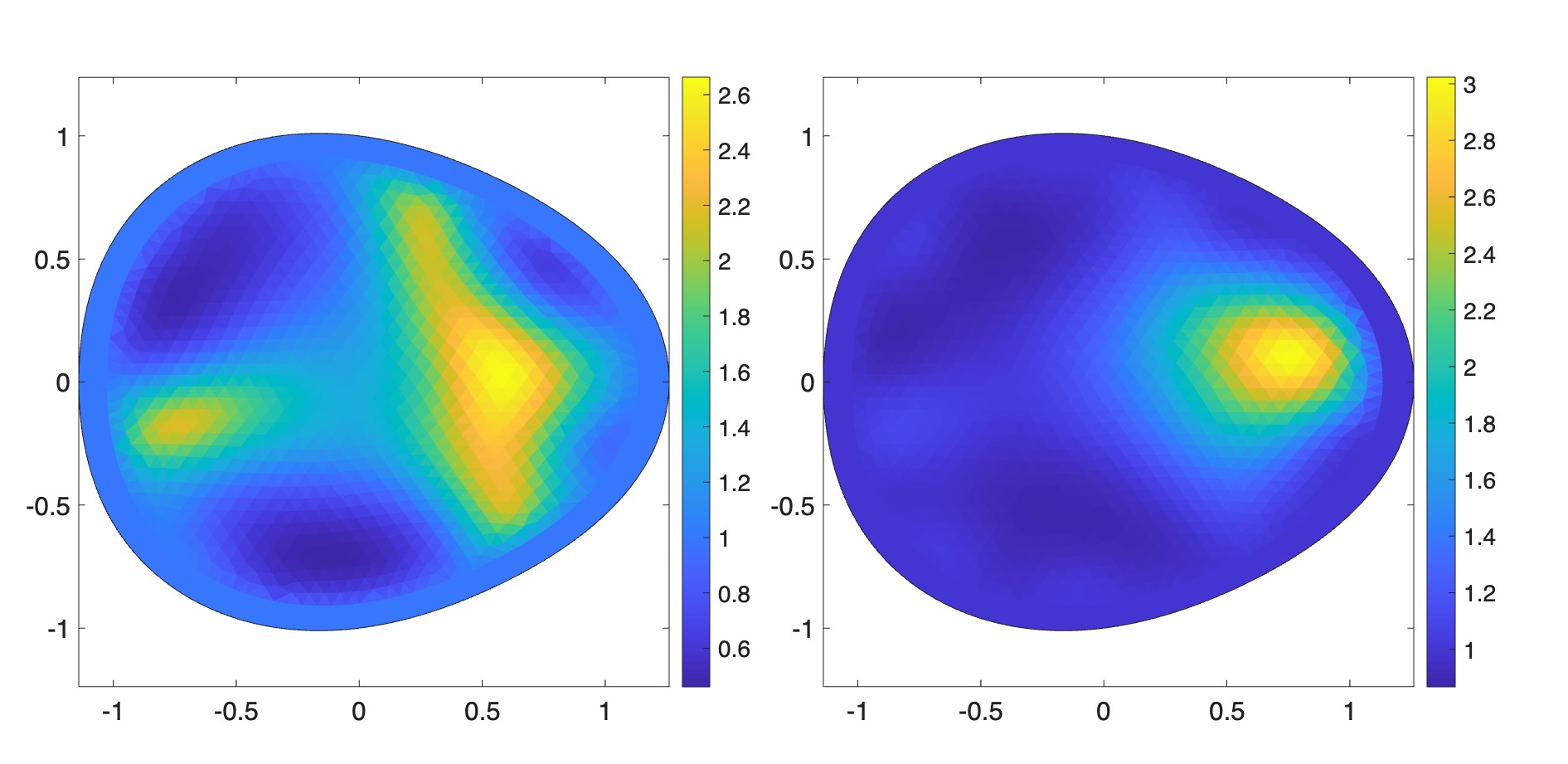}
}
\caption{On the left, reconstruction obtained by ignoring the approximation error. The geometry artifacts are strongly present in this reconstructions. In the right panel, the reconstruction is computed by using a the spotlight projection calculated with only five random draws of the approximation error.}\label{fig:EIT reconstructions}
\end{figure}

To illustrate the performance of the proposed approaches to model reduction, highlighting the importance of accounting for the approximation errors, we use the Gauss-Newton iteration in the same manner on the unprojected problem, effectively neglecting to account for the geometric uncertainty.  In the regularized solution of the inverse problem neglecting to account for the approximation error, displayed in the left panel of Figure~\ref{fig:EIT reconstructions}, the artifacts due to the incorrect geometry are clearly visible.

To mitigate the approximation error artifacts, we estimated the approximation error based on a sample of only five random realizations, and proceeded to perform the spotlight projection outlined above using the same Gauss-Newton approach. The reconstruction, shown on the right panel of Figure~\ref{fig:EIT reconstructions}, free of artifacts, also matches quite accurately the dynamic range.  

\section{Discussion}

This article is bridging two different approaches to deal with model uncertainties, the Bayesian approximation error paradigm and the clutter removal by spotlight projections. The two approaches have some fundamental differences, the former building on Bayesian analysis of inverse problems, while the latter being purely linear algebraic. The approach adopted in this paper is navigating between these two approaches, taking advantage of the Bayesian paradigm in building samples of the approximation error, and the linear algebraic approach by basing the inverse solution purely on projections and on solving regularized least squares problems. The idea of using the prior probability densities to generate samples of the approximation error leads naturally to a formalism that has a close similarity to the sketching methods in randomized linear algebra, and gives rise to a novel way of seeing sketching in the context of Bayesian inverse problems as priorsketching. Further connections, e.g., with approaches like Hessian sketching to solve inverse problems will be the topic of future work. To keep the discussion simple, the analysis of the BAE was limited to the enhanced error model formalism, in which the interdependency of the approximation error and the quantity of interest is ignored. Discussion of the spotlight approach in the light of more sophisticated models of BAE is also a topic of interest for future studies.

{\bf{Acknowledgments:}}
This material is based upon work supported by the National Science Foundation under Grant No. DMS-1929284 while the authors were in residence at the Institute for Computational and Experimental Research in Mathematics (ICERM) in Providence, RI, during the "Stochastic and Randomized Algorithms in Scientific Computing: Foundations and Applications" semester program. 
The work of DC was partly supported by the NSF grant DMS 2513481, as well as by the Simons Foundation Fellowship, and that of  ES by the NSF grants DMS 2204618 and DMS 2513481. The support is acknowledged with gratitude.

\bibliographystyle{siam} 
\bibliography{Biblio}

@article{arridge2006approximation,
  title={Approximation errors and model reduction with an application in optical diffusion tomography},
  author={Arridge, Simon R and Kaipio, Jari P and Kolehmainen, Ville and Schweiger, Martin and Somersalo, Erkki and Tarvainen, Tanja and Vauhkonen, Marko},
  journal={Inverse problems},
  volume={22},
  number={1},
  pages={175--195},
  year={2006}
}

@article{nissinen2010compensation,
  title={Compensation of modelling errors due to unknown domain boundary in electrical impedance tomography},
  author={Nissinen, Antti and Kolehmainen, Ville Petteri and Kaipio, Jari P},
  journal={IEEE transactions on medical imaging},
  volume={30},
  number={2},
  pages={231--242},
  year={2010},
  publisher={IEEE}
}

@article{koulouri2016compensation,
  title={Compensation of domain modelling errors in the inverse source problem of the Poisson equation: Application in electroencephalographic imaging},
  author={Koulouri, Alexandra and Rimpil{\"a}inen, Ville and Brookes, Mike and Kaipio, Jari P},
  journal={Applied Numerical Mathematics},
  volume={106},
  pages={24--36},
  year={2016},
  publisher={Elsevier}
}

@article{huttunen2007approximation,
  title={Approximation errors in nonstationary inverse problems},
  author={Huttunen, Janne MJ and Kaipio, Jari P and Somersalo, E},
  journal={Inverse problems and imaging},
  volume={1},
  number={1},
  pages={77},
  year={2007},
  publisher={American Institute of Mathematical Sciences}
}

@article{kaipio2019bayesian,
  title={A Bayesian approach to improving the Born approximation for inverse scattering with high-contrast materials},
  author={Kaipio, Jari P and Huttunen, Tomi and Luostari, Teemu and L{\"a}hivaara, Timo and Monk, Peter B},
  journal={Inverse Problems},
  volume={35},
  number={8},
  pages={084001},
  year={2019},
  publisher={IOP Publishing}
}

@dataset{bubbaOA,
  author       = {Bubba, Tatiana A. and
                  Hauptmann, Andreas and
                  Huotari, Simo and
                  Rimpeläinen, Juho and
                  Siltanen, Samuli},
  title        = {Tomographic {X}-ray data of a lotus root filled with
                   attenuating objects
                  },
  month        = sep,
  year         = 2016,
  publisher    = {Zenodo},
  version      = {1.0.0},
  doi          = {10.5281/zenodo.1254204},
  url          = {https://doi.org/10.5281/zenodo.1254204},
}

@book{kaipio2005statistical,
  title={Statistical and computational inverse problems},
  author={Kaipio, Jari P and Somersalo, Erkki},
  year={2005},
  publisher={Springer}
}

@article{kaipio2007statistical,
  title={Statistical inverse problems: discretization, model reduction and inverse crimes},
  author={Kaipio, Jari and Somersalo, Erkki},
  journal={Journal of computational and applied mathematics},
  volume={198},
  number={2},
  pages={493--504},
  year={2007},
  publisher={Elsevier}
}

@book{calvetti2023bayesian,
  title={Bayesian scientific computing},
  author={Calvetti, Daniela and Somersalo, Erkki},
  year={2023},
  publisher={Springer}
}

@article{jaaskelainen2024projection,
author = {J\"{a}\"{a}skel\"{a}inen, A. and Toivanen, J. and H\"{a}nninen, A. and Kolehmainen, V. and Hyv\"{o}nen, N.},
title = {Projection-Based Preprocessing for Electrical Impedance Tomography to Reduce the Effect of Electrode Contacts},
journal = {SIAM Journal on Imaging Sciences},
volume = {18},
number = {3},
pages = {1681-1706},
year = {2025},
doi = {10.1137/24M1719517},

URL = { 
    
        https://doi.org/10.1137/24M1719517
    
    

},
eprint = { 
    
        https://doi.org/10.1137/24M1719517
    
    

}
}

@article{somersalo1992existence,
  title={Existence and uniqueness for electrode models for electric current computed tomography},
  author={Somersalo, Erkki and Cheney, Margaret and Isaacson, David},
  journal={SIAM Journal on Applied Mathematics},
  volume={52},
  number={4},
  pages={1023--1040},
  year={1992},
  publisher={SIAM}
}

@article{calvetti2025spotlight,
  title={Spotlight inversion by orthogonal projections},
  author={Calvetti, Daniela and Hyv{\"o}nen, Nuutti and Kolehmainen, Ville and Somersalo, Erkki},
  journal={arXiv preprint arXiv:2509.15512},
  year={2025}
}

@article{calvetti2025sparsity,
  title={Sparsity-promoting hierarchical Bayesian model for EIT with a blocky target},
  author={Calvetti, Daniela and Pragliola, Monica and Somersalo, Erkki},
  journal={Journal of Computational Physics},
  pages={114255},
  year={2025},
  publisher={Elsevier}
}

@article{kolehmainen2013recovering,
  title={Recovering boundary shape and conductivity in electrical impedance tomography},
  author={Kolehmainen, Ville and Lassas, Matti and Ola, Petri and Siltanen, Samuli},
  journal={Inverse Problems and Imaging},
  volume={7},
  number={1},
  pages={217--242},
  year={2013},
  publisher={American Institute of Mathematical Sciences}
}

@article{kolehmainen2008electrical,
  title={Electrical impedance tomography problem with inaccurately known boundary and contact impedances},
  author={Kolehmainen, Ville and Lassas, Matti and Ola, Petri},
  journal={IEEE transactions on medical imaging},
  volume={27},
  number={10},
  pages={1404--1414},
  year={2008},
  publisher={IEEE}
}

@article{hyvonen2017compensation,
  title={Compensation for geometric modeling errors by positioning of electrodes in electrical impedance tomography},
  author={Hyv{\"o}nen, N and Majander, H and Staboulis, Stratos},
  journal={Inverse Problems},
  volume={33},
  number={3},
  pages={035006},
  year={2017},
  publisher={IOP Publishing}
}

@article{candiani2019computational,
  title={Computational framework for applying electrical impedance tomography to head imaging},
  author={Candiani, Valentina and Hannukainen, Antti and Hyvonen, Nuutti},
  journal={SIAM Journal on Scientific Computing},
  volume={41},
  number={5},
  pages={B1034--B1060},
  year={2019},
  publisher={SIAM}
}

@article{calvetti2008hypermodels,
  title={Hypermodels in the Bayesian imaging framework},
  author={Calvetti, Daniela and Somersalo, Erkki},
  journal={Inverse Problems},
  volume={24},
  number={3},
  pages={034013},
  year={2008}
}

@article{calvetti2018iterative,
  title={Iterative updating of model error for Bayesian inversion},
  author={Calvetti, Daniela and Dunlop, Matthew and Somersalo, Erkki and Stuart, Andrew},
  journal={Inverse Problems},
  volume={34},
  number={2},
  pages={025008},
  year={2018},
  publisher={IOP Publishing}
}

@article{kilmer2001choosing,
  title={Choosing regularization parameters in iterative methods for ill-posed problems},
  author={Kilmer, Misha E and O'Leary, Dianne P},
  journal={SIAM Journal on matrix analysis and applications},
  volume={22},
  number={4},
  pages={1204--1221},
  year={2001},
  publisher={SIAM}
}

@article{kaipio2013approximate,
  title={Approximate marginalization over modeling errors and uncertainties in inverse problems},
  author={Kaipio, Jari and Kolehmainen, Ville},
  journal={Bayesian theory and applications},
  pages={644--672},
  year={2013},
  publisher={Oxford University Press Oxford}
}

@article{nissinen2011reconstruction,
  title={Reconstruction of domain boundary and conductivity in electrical impedance tomography using the approximation error approach},
  author={Nissinen, Antti and Kolehmainen, Ville and Kaipio, Jari P},
  journal={Visualization of Mechanical Processes: An International Online Journal},
  volume={1},
  number={3},
  year={2011},
  publisher={Begel House Inc.}
}

\end{document}